\documentclass[11pt,reqno]{amsart}
\usepackage{amssymb}
\usepackage{amsfonts}
\usepackage{amsthm}
\numberwithin{equation}{section}
\theoremstyle{plain}
\newtheorem{theorem}{Theorem}[section]
\newtheorem{proposition}[theorem]{Proposition}

\theoremstyle{definition}
\newtheorem{definition}[theorem]{Definition}

\newcommand{\refE}[1]{(\ref{E:#1})}
\newcommand{\refS}[1]{Section~\ref{S:#1}}

\newcommand{\C}{\ensuremath{\mathbb{C}}}
\newcommand{\N}{\ensuremath{\mathbb{N}}}

\renewcommand{\P}{\ensuremath{\mathbb{P}}}
\newcommand{\Z}{\ensuremath{\mathbb{Z}}}
\newcommand{\K}{\ensuremath{\mathbb{K}}}

\newcommand{\cins}{\frac 1{2\pi\mathrm{i}}\int_{C_S}}

\newcommand{\cintl}[1]{\frac 1{24\pi\mathrm{i}}\int_{#1 }}
\newcommand{\g}{\mathfrak{g}}
\newcommand{\gb}{\overline{\mathfrak{g}}}
\newcommand{\gh}{\widehat{\mathfrak{g}}}
\newcommand{\Gb}{\overline{\mathcal{G}}}

\newcommand{\A}{\mathcal{A}}
\newcommand{\W}{\mathcal{W}}
\newcommand{\V}{\mathcal{V}}

\renewcommand{\L}{\mathcal{L}}

\renewcommand{\H}{\mathrm{H}}

\newcommand{\eee}{(e_1-e_2)(e_1-e_3)}
\newcommand{\ddX}{\frac {d}{dX}} 
 
\newcommand{\ddz}{\frac {d}{dz}} 
\newcommand{\tensor}{\otimes}   
\newcommand{\htensor}{\widehat{\otimes}}

\begin{document}
\vspace*{-1cm}
\hbox{ }
{{\hspace*{\fill} math/0610851}}%

\vspace*{2cm}

\title[Global Geometric Deformations \quad, \today] 
{Global Geometric Deformations of the Virasoro algebra, current and 
affine  algebras by
Krichever-Novikov type algebras}
\thanks{}
\author[A. Fialowski]{Alice Fialowski}
\address[Alice Fialowski]{Institute of Mathematics,
E\"otv\"os Lor\'and University, P\'azm\'any P\'eter s\'et\'any 1/C,
H-1117 Budapest, Hungary}
\email{fialowsk@cs.elte.hu}
\author[M. Schlichenmaier,\quad \today ]{Martin Schlichenmaier}
\address[Martin Schlichenmaier]{University of Luxembourg,
Institute of Mathematics, Campus Limpertsberg, 
162 A, Avenue de la Faiencerie,
L-1511 Luxembourg, Grand-Duchy of Luxembourg}
\email{Martin.Schlichenmaier@uni.lu}
\begin{abstract}
In two earlier articles we 
constructed algebraic-geometric families of genus one (i.e. elliptic)
Lie algebras of Krichever-Novikov type. The considered algebras are 
vector fields, current and affine Lie algebras.
These families deform the Witt algebra, 
the Virasoro algebra, the classical current, and the 
affine Kac-Moody Lie algebras respectively.
The constructed families are 
not equivalent (not even locally) to the trivial families, despite the fact that
the classical algebras are formally rigid.
This effect is due to the fact that the algebras
are infinite dimensional.
In this article the results are reviewed and developed 
further.
The constructions are induced by the geometric process of
degenerating the elliptic curves to singular cubics.
The algebras  are of relevance in the global operator
approach to the Wess-Zumino-Witten-Novikov models 
appearing in the quantization of Conformal Field Theory.
%
%
\end{abstract}
\date{27.10.2006}
\maketitle

\vskip 1.0cm
\section{Introduction}\label{S:intro}

Deformation theory plays a crucial role in all branches of mathematics and
physics. 
In physics the mathematical theory of deformations has proved to be a powerful
tool in modeling physical reality. 
The concepts {\it symmetry} and {\it deformations} are considered
to be two fundamental guiding principle for developing the
physical theory further.
From the mathematical point of view
considering deformed objects will give additional information about the 
original object itself. In particular, how is its relation to
``neighbouring'' objects. This can be made precise  with the
notion of moduli space, classifying inequivalent objects of the same
type. The moduli space should be equipped with a ``geometric''
structure such that ``nearby points''  should be
also ``nearby'' in the sense of deforming the structure of the
initial object.
Moreover, assuming that such a moduli space exists, its dimension
should be equal to the number of inequivalent  deformation directions.
Maybe there exists even a deformation family containing 
every possible deformation.

Clearly, this general remarks are rather vague. To make them more 
precise first one has to be more precise about the structure to deform.
A very famous and well-developed domain is the deformation theory
of  complex analytic structures of a compact complex manifold $M$.
We do not have the place to recall here this theory, but refer only to 
\cite{Palm} for results and details.
Let us only mention that a fundamental role is played by the
first cohomology space $\H^1(M,T_M)$ of $M$ with values in the
holomorphic tangent sheaf $T_M$.
In particular, if this space is trivial, $M$ will be rigid, i.e. it cannot
be deformed in something which is not isomorphic to it.

Here we will deal with 
deformations of Lie algebras, in particular of such of infinite dimension.
 Formal deformations of arbitrary rings and
 associative algebras, and the related cohomology questions,
 were first investigated by Gerstenhaber, in a series of articles
 \cite{Ger}. The notion of deformation was applied to
 Lie algebras by Nijenhuis
 and Richardson \cite{NuijRich1}, \cite{NuijRich2}.

The cohomology space related to 
the deformations of a Lie algebra $\L$ is the Lie algebra
two-cohomology $\H^2(\L,\L)$ of $\L$ with  
values in the adjoint module.
We will explain this in \refS{defcom}.
As long as the Lie algebra 
is finite-dimensional, the relation is rather tight. In particular, if
the cohomology space vanishes, the Lie algebra will be rigid in 
all respects.

But the algebras which  are e.g. of relevance in Conformal 
Field Theory, integrable systems related to partial differential
equations, etc. are typically infinite dimensional.  We are interested 
here in these algebras.
We showed in two articles \cite{FiaSchlwitt},  \cite{FiaSchlaff}
that the relation to cohomology is not so tight anymore.
In particular we constructed nontrivial geometric 
deformation families for the Witt algebra (resp. its universal
central extension the Virasoro algebra) and for the
current algebras (resp. their central extensions the affine algebras),
despite the fact that 
the cohomology spaces for those algebras are 
trivial and hence the algebras are formally rigid
\cite{Fia2}, \cite{LeRog}.
This is a phenomena which in finite dimension cannot occur.

\medskip

Here we report on our results and the constructions to be found in 
 \cite{FiaSchlwitt},  \cite{FiaSchlaff} 
and continue our investigation.
The Witt algebra is the algebra consisting of those 
meromorphic vector fields on the Riemann sphere which are holomorphic
outside $\{0,\infty\}$. A basis and the associated structure is 
given by
\begin{equation*}
l_n=z^{n+1}\ddz,\ n\in\Z,
\quad
\text{ with Lie bracket}
\quad 
[l_n,l_m]=(m-n)\,l_{n+m}.
\end{equation*}
The Virasoro algebra is its universal central extension
\begin{equation*}
[l_n,l_m] = (m-n)l_{n+m} + \tfrac{1}{12}(m^3-m)\delta_{n,-m}\,t,
\qquad
[l_n,  t] =0,
\end{equation*}
with $t$ an additional basis element which is central.

Furthermore we consider the case of 
{current algebras}
$\gb=\g\otimes \C[z^{-1},z]$ and their central extensions
$\gh$, the {\it affine Lie algebras}.
Here $\g$ is a finite-dimensional Lie algebra (which for simplicity
we assume to be simple).
With the Cartan-Killing  form $\beta$
the central extension $\gh$ is the vector space 
$\gb\oplus t\,\C$ endowed with the Lie bracket
\begin{equation*}
[x\otimes z^n, y\otimes z^m]=[x,y]\otimes z^{n+m}-\beta(x,y)\cdot 
n\cdot\delta_{m}^{-n}\cdot t,
\quad [t,\gh]=0,
\quad x,y\in\g,\ n,m\in\Z.
\end{equation*}
As already mentioned,
these algebras are 
rigid.

The families we construct,  appear as families of higher-genus
multi-point algebras of Krichever-Novikov type,
see \refS{kn} for their definitions.
Hence, they are not just abstract families, but families
obtained by geometric processes. 
The results obtained do not have only relevance in
deformation theory of algebras, but they also are of importance in
areas where vector fields, current, and affine algebras play a role.

\medskip 
A very prominent application
is two-dimensional conformal field theory (CFT) and its quantization.
It is well-known that the Witt algebra, 
the Virasoro algebra, the current algebras, the 
affine algebras, and their representations 
are of fundamental importance  for CFT
on the Riemann sphere
(i.e. for genus zero), see \cite{BP}.
Krichever and Novikov \cite{KNFac} proposed in  the case of
higher genus Riemann surfaces (with two insertion 
points) the use of global operator fields 
which are given with the help of  the Lie algebra 
of vector fields of Krichever-Novikov type,
certain related algebras, 
and their representations (see \refS{kn} below).

Their approach was extended by 
Schlichenmaier to the multi-point
situation (i.e. an arbitrary number of insertion points was allowed)
\cite{SchlDiss}, \cite{Schlkn}, \cite{Schleg}, \cite{Schlce}.
The necessary central extensions where constructed. Higher
genus multi-point current and affine algebras were introduced
\cite{SchlCt}.
These algebras consist of meromorphic objects on a Riemann surface
which are holomorphic outside a finite set $A$  of points.
The set $A$ is divided into two disjoint subsets
$I$ and $O$. With respect to some possible interpretation of the
Riemann surface as the world-sheet of a string, the points
in $I$ are called {\it in-points}, the points in $O$ 
are called {\it out-points}, 
corresponding to incoming and outgoing free strings.
The world-sheet itself corresponds to possible interaction.
This splitting introduces an almost-graded structure (see \refS{kn})
for the algebras and their representations.
Such an almost-graded structure is needed to construct representations
of relevance in the context of the quantization of CFT, e.g. highest weight
representations, 
fermionic Fock space representations, etc.

In the following we give  more information on 
a special model.
In the process of quantization of  conformal fields one has
to consider families of algebras and representations over the
moduli space of compact Riemann surfaces
(or equivalently, of  smooth projective curves over $\C$) of genus $g$ with
$N$ marked points.
Models of most importance in CFT are the 
Wess-Zumino-Witten-Novikov models (WZWN).
Tsuchiya, Ueno, and Yamada \cite{TUY} gave a sheaf version of
WZWN models over the moduli space.
In \cite{SchlShwznw}, 
\cite{SchlSh}
Schlichenmaier and Sheinman  developed a global operator version.
In this context
of particular interest is the situation
$I=\{P_1,\ldots, P_K\}$, the marked points we want to vary,
 and $O=\{P_\infty\}$, a reference point. We obtain
families of algebras over 
the moduli space $\mathcal{M}_{g,K+1}$ of curves of genus $g$ with
$K+1$ marked points, and we are exactly in the middle of the
main subject of this article.
In \cite{SchlShwznw} and   \cite{SchlSh} it is shown that
there exists a global operator description of 
WZWN  models with the help of the Krichever
Novikov objects at least over  a dense open subset of the moduli space.
Starting from families of representations $\V$ of families of
higher genus affine algebras 
(see \refS{kn} for their definition) 
the vector bundle  of conformal blocks can 
be defined as the vector bundle with fibre 
(over the moduli point $b=[(M,\{P_1,\ldots, P_K \},\{P_\infty\})]$)
the quotient space of the fibre $\V_b$ modulo the subspace generated by
the vectors obtained by the action
 of those elements of the affine algebra  which vanish at the
reference point $P_\infty$ (i.e. the fibre is the space of 
coinvariants of this subalgebra).

The bundle of conformal blocks carries a connection
called
the {\it Kniz\-hnik-Zamolod\-chikov connection}.
In its definition 
an important role is played by the Sugawara construction which 
associates to representations of affine algebras representations
of the (almost-graded) centrally extended vector
field algebras, see \cite{SchlShSug}.
A certain subspace of the vector field algebra  (assigned
to the moduli point $b$) corresponds to tangent directions
on the moduli space $\mathcal{M}_{g,K+1}$ at the point $b$.

Now clearly, the following question is of fundamental importance. What happens
if we approach the boundary of the moduli space?
The boundary components correspond to curves with singularities.
Resolving the singularities yields curves of lower genera.
By geometric degeneration we obtain families of (Lie) algebras containing
a lower genus algebra (or sometimes a subalgebra of it), 
corresponding to a suitable collection of marked points,
as special element.
Or reverting the perspective,
we obtain a typical situation of the deformation  of an algebra
corresponding in some way to a lower genus situation, containing
higher genus algebras as the other elements in the family.
Such kind of geometric degenerations
are of fundamental importance if one wants to prove Verlinde
type formula via factorization and normalization technique,
see \cite{TUY}.

By a maximal degeneration a collection of $\P^1(\C)$'s will 
appear.
Indeed, the examples considered in this article are exactly of this type.
The deformations  appear as families 
of vector fields and current algebras which are naturally defined over the
moduli space 
of genus one curves (i.e. of elliptic curves, or
equivalently of complex one-dimensional tori) with two marked points.
These deformations are
associated to  geometric degenerations of elliptic curves to
singular cubic curves. 
The desingularization (or normalization) of their singularities will yield 
the projective line as normalization. We will end up with
algebras related to the genus zero case.
The full geometric picture behind the degeneration 
 was discussed in \cite{Schltori}.
In particular, we like to point out, that even if one starts with
two marked points, by passing  to the boundary of the moduli space
one is forced to consider more points (now for a curve of lower
genus).

\bigskip
{\it Acknowledgements:}
The authors thank 
the {\it Mathematisches Forschungsinstitut in Oberwolfach (MFO)} for
the opportunity to organize together with Marc de Montigny and Sergey
P. Novikov a workshop on
``Deformations and Contractions in Mathematics and Physics''
in January 2006.
The workshop brought together mathematicians and physicists at a 
very inspiring meeting.
Furthermore, we like to thank the 
Erwin-Schr\"odinger Institute (ESI) in Vienna, 
 the
Institut des Hautes \'Etudes Scientifiques (IHES)
in Bures-sur-Yvette were certain parts of the work was carried out.
The work was partially supported by
grants OTKA T043641 and T043034 and the research grant
R1F105L15 of the University of Luxembourg.


\section{Deformations of Lie algebras}
\label{S:def}
In the physics literature a Lie algebra $\L$ is often given in
terms of generators and structure constants.
Let $V$ be a finite- or infinite dimensional 
complex vector space  with basis
$\{T_a\}_{a\in J}$ then a Lie algebra structure on $V$ 
can be given by the  structure equations, i.e. the Lie bracket,
\begin{equation}\label{E:lieb}
[T_a,T_b]=\sum_{c\in J}{}'\; C_{a,b}^c T_c,\qquad  a,b\in J,
\end{equation}
with structure constants  $C_{a,b}^c\in\C$.
The symbol  $\sum'$ denotes that for fixed $a,b\in J$ only
for finitely many $c$  the coefficient $C_{a,b}^c\ne 0$.
In terms of structure constants the 
necessary and sufficient conditions for $\L$ being a 
Lie algebra (i.e. the anti-symmetry and the Jacobi identity) can be 
written as
\begin{equation}\label{E:structlie}
\begin{gathered}
C_{a,b}^c+C_{b,a}^c=0,\qquad  a,b,c\in J,
\\
\sum_{l\in J}
\left(
C_{a,b}^l C_{l,c}^d+
C_{b,c}^l C_{l,a}^d+
C_{c,a}^l C_{l,b}^d
\right)=0,
\qquad  a,b,c,d\in J.
\end{gathered}
\end{equation}
Deforming the Lie algebra structure corresponds intuitively to
making the system of coefficients $\{C_{a,b}^c\}$ depending on one or more
parameters.

In a more compact manner a  Lie algebra $\L$, 
i.e. its bracket $[.,.]$, might be written 
with an anti-symmetric bilinear form
$$\mu_0:\L\times \L\to \L, \qquad \mu_0(x,y)=[x,y],$$
fulfilling certain additional conditions corresponding 
to  the Jacobi identity.
Consider on the same vector space $\L$ is modeled on,  a family
of Lie structures
\begin{equation}\label{E:naivdef}
\mu_t=\mu_0+t\cdot \phi_1+t^2\cdot \phi_2+\cdots  \ ,
\end{equation}
with bilinear maps $\phi_i:\L\times \L \to \L$ such that 
$\L_t:=(\L,\mu_t)$ is a Lie algebra and $\L_0$ is the
Lie algebra we started with.
The family $\{\L_t\}$ is a \emph{deformation} of $\L_0$.

Up to this point  we did not specify the ``parameter'' $t$.
Indeed,  different choices are possible.
\begin{enumerate}
\item
The parameter
$t$  might be a  variable which allows to plug in  
numbers $\alpha\in\C$. 
In this case $\L_{\alpha}$ is a Lie algebra for every $\alpha$ for which
the expression \refE{naivdef} is defined.
The family can be considered as deformation over the affine line
$\C[t]$ or over the convergent power series $\C\{\!\{t\}\!\}$.
The deformation is called a {\it geometric} or an 
{\it analytic deformation}
respectively.
\item
We consider $t$ as a formal variable and we allow 
infinitely many terms in \refE{naivdef}. It might be the case that
$\mu_t$ does not exist if we plug in for $t$ any other value  different from 
$0$. 
In this way  we obtain deformations over the ring of formal 
power series $\C[[t]]$. The corresponding deformation is a
{\it formal deformation}.
\item
The parameter
$t$ is considered as an infinitesimal variable, i.e. we take 
$t^2=0$. We obtain {\it infinitesimal deformations} defined 
over the quotient $\C[X]/(X^2)=\C[[X]]/(X^2)$.
\end{enumerate}
We could even consider more general situations for the
parameter space.
See  Appendix A for a  general mathematical
definition of a deformation.

\medskip

There is always the trivially deformed family given by
$\mu_t=\mu_0$ for all values of $t$.
Two families $\mu_t$ and 
$\mu_t'$ deforming the same $\mu_0$
are \emph{equivalent} if there exists a linear
automorphism
(with the same vagueness about the meaning of $t$)
\begin{equation}
\psi_t=id+t\cdot \alpha_1+t^2\cdot \alpha_2 +\cdots
\end{equation}
with $\alpha_i:\L\to\L$ linear maps such that 
\begin{equation}
\mu_t'(x,y)=\psi_t^{-1}(\mu_t(\psi_t(x),\psi_t(y))).
\end{equation}
A Lie algebra $(\L,\mu_0)$ is called \emph{rigid} if every deformation
$\mu_t$ of $\mu_0$ is locally equivalent to the trivial family.
Intuitively, this says that $\L$ cannot be deformed.

The word
``locally'' in the definition of rigidity means 
that we only consider the situation 
for $t$ ``near 0''.
Of course, this depends on the category we consider.
As on the formal and  the infinitesimal level there exists only
one closed point, i.e. the point $0$ itself, every deformation
over $\C[[t]]$ or $\C[X]/(X^2)$ is already local.
This is different on the geometric and analytic level.
Here it means that there exists an open 
neighbourhood $U$ of $0$ such that the family restricted to it is 
equivalent to the trivial one.
In particular, this implies $\L_\alpha\cong \L_0$ for all $\alpha\in U$.

Clearly, a question of fundamental interest is to decide whether a 
given Lie algebra is rigid. Moreover, the question of rigidity will  depend
on the category we consider. Depending on the set-up we will
have to consider infinitesimal, formal, geometric, and analytic 
rigidity.
If the
 algebra is not rigid, one would like to know whether there exists
a moduli space of (inequivalent) deformations. If so, what is
its structure, dimension, etc?

As explained in the introduction, deformation problems and moduli space
problems are related to adapted cohomology theories.
To a certain extend (in particular for the finite-dimensional case)
this is also true for deformations of Lie algebras.
But as far as geometric and algebraic deformations are concerned
it is wrong for infinite dimensional Lie algebras as our 
examples show.

\medskip
The corresponding relations to cohomology will be explained in 
\refS{defcom}. 
To see later why the results are  different in the infinite dimensional
case, let us first discuss the finite-dimensional case.
Let $\L$ be a finite-dimensional Lie algebra of
dimension $n$ over $\C$ and denote the underlying vector space by $V$.
The structure constants
from \refE{lieb}
$\{C_{a,b}^c\}_{a,b,c=1,\ldots, n}$ are elements of $\C^{n^3}$.
The conditions 
\refE{structlie}, 
which are necessary and sufficient that \refE{lieb} defines  a Lie algebra,
are algebraic 
equations. 
The vanishing set $Lalg_n$ (i.e. the set of common zeros of these
equations)
in $\C^{n^3}$ ``parameterizes'' the possible Lie algebra structures on
the $n$-dimensional vector space $V$.
As the conditions are algebraic the vanishing set will be
a (not necessarily irreducible) variety. In fact
it would be better to talk about  $Lalg_n$ as a
scheme, as one should better consider the not necessarily
reduced structure on  $Lalg_n$.

The Lie structure 
$\mu$ is a bilinear map
$V\times V\to V$  and 
the structure constants might be considered 
as elements of $V^*\otimes V^*\otimes V$ with
 $V^*$ the dual space of $V$.

If  we make a change of basis, the structure constants will change.
The two set of structure constants will define 
isomorphic Lie algebras.
The corresponding effect can be described by 
a linear automorphism  $\Phi\in \mathrm{Gl}(V)$.
It will define 
an action   on  $V^*\otimes V^*\otimes V$ by
\begin{equation}
(\Phi\star\mu)(x,y)=\Phi(\mu(\Phi^{-1}(x),\Phi^{-1}(y))).
\end{equation}
If $\mu$ corresponds to a Lie algebra structure, 
$\Phi\star\mu$ will also be a Lie algebra. Hence
$\Phi\star$ will be an action on $Lalg_n$

The Lie algebras $(V,\mu)$ and  $(V,\mu')$
are isomorphic 
iff  $\mu$ and $\mu'$ are  
in the same orbit under this  $\mathrm{Gl}(V)$ action.
On the level of structure constants, i.e. after fixing a
basis in $V$, we obtain a $\mathrm{Gl}(n)$ action on
 $Lalg_n$.
In this way the isomorphy classes of Lie algebras of dimension
$n$ correspond exactly to the  $\mathrm{Gl}(n)$ orbits of $Lalg_n$.

The variety $Lalg_n$  decomposes  into different orbits under
the $\mathrm{Gl}(n)$-action. 
Let $x_0$ be a point in $Lalg_n$ (defining the Lie structure $\mu_0$).
All ``nearby'' Lie structures $\mu$ correspond to points $x$ near $x_0$.
Of course the $\mathrm{Gl}(n)$ orbit of $x_0$ passes through $x_0$.
If all points in an open neighbourhood of $x_0$ lie in this orbit then
this implies that all ``nearby'' Lie structures $\mu_t$ are 
isomorphic to $\mu_0$.
A Lie algebra is called rigid in the orbit sense, if the corresponding orbit 
is Zariski open in $Lalg_n$%
\footnote{A subset is called Zariski open if it is the complement of 
the vanishing set of finitely many algebraic equations. Zariski open 
subset are always open in the usual topology.}.
In particular, rigidity in the orbit sense implies rigidity in the
geometric and analytic sense.

Intuitively the ``moduli space'' of finite-dimensional Lie structures should
correspond to the orbit space under the $\mathrm{Gl}(n)$-action. 
But as in the boundary  of certain
orbits there might be  different orbits this will need some
modification. Indeed, the 
problem of the geometric structure  of 
the ``moduli space'' is  rather delicate and as we will not
need it here,  we will 
not discuss it, see Bjar and Laudal \cite{BL}.

\bigskip

Back to our Lie algebras of arbitrary dimensions.
Special types of deformations are {\it jump deformations}.
They are typically given as families over  a parameters space 
(parameterized e.g. by $t$) 
around $0$, such that $\L_t\cong \L_{t'}$ 
as long as $t,t'\ne 0$,
but 
 $\L_t\not\cong \L_{0}$. In the finite-dimensional case,
 considered above,
the element $\L_0$ 
will be necessarily 
a boundary point of the orbit of  $\L_t$, $t\ne 0$ which does not lie 
in the orbit itself.
This says it is an element in the Zariski closure 
of the orbit but not of the orbit itself.
Sometimes in physics one talks about {\it contractions}. This language is dual 
to the language of jump deformations. Here $\L_0$ is a contraction of 
the isomorphy type of  $\L_t$, $t\ne 0$. 
Moreover, in the finite-dimensional case the 
possible contractions of $\L$ are given by the boundary points  of its
$\mathrm{Gl}(n)$ orbit.


\section{Cohomological description}
\label{S:defcom}
For Lie algebra deformations the relevant cohomology space is
$\H^2(\L,\L)$, the space of Lie algebra two-cohomology classes with
values in the adjoint module $\L$.

Recall that these cohomology classes are classes of two-cocycles
modulo coboundaries.
An antisymmetric bilinear 
map $\phi:\L\times \L\to\L$ is a Lie algebra \emph{two-cocycle} if
$d_2\phi=0$, or expressed explicitely
\begin{equation}
\phi([x,y],z)+
\phi([y,z],x)+
\phi([z,x],y)-
[x,\phi(y,z)]+
[y,\phi(z,x)]-
[z,\phi(x,y)]=0.
\end{equation}
The map $\phi$ will be 
a \emph{coboundary} if there exists a linear map $\psi:\L\to\L$ with
\begin{equation}
\phi(x,y)=(d_1\psi)(x,y):=
\psi([x,y])
-[x,\psi(y)]
+[y,\psi(x)].
\end{equation}
If we write down  the Jacobi identity for $\mu_t$  
given by \refE{naivdef}
then it can be immediately verified 
that the 
first non-vanishing $\phi_{i}$ has to be  a two-cocycle
in the above sense.
Furthermore, if $\mu_t$ and $\mu_t'$ are equivalent then
the corresponding 
 $\phi_{i}$  and  $\phi_{i}'$ are cohomologous, i.e. their
difference is a coboundary.

The following results are well-known:
\begin{enumerate}
\item
$\H^2(\L,\L)$ classifies infinitesimal deformations of $\L$ \cite{Ger}.
\label{it1}
\item
If $\dim \H^2(\L,\L)<\infty$ then all formal
deformations of $\L$ up to equivalence can be realised in this vector space
\cite{FiaFuMini}.\label{it2}
\item
If $\H^2(\L,\L)=0$ then $\L$ is infinitesimally and formally rigid
(this follows directly from (\ref{it1})  and (\ref{it2})).
\item
If $\dim \L<\infty$ then $\H^2(\L,\L)=0$ implies 
that $\L$ is also rigid in the geometric and analytic  sense
\cite{Ger}, \cite{NuijRich1}, \cite{NuijRich2}.
\label{it4}
\end{enumerate}

As our examples show,
without the  condition $\dim \L<\infty$
point \ref{it4} is not true anymore.

For the Witt algebra $\W$ one has $\H^2(\W,\W)=0$ 
(\cite{Fia2}, see also \cite{FiaSchlwitt}).
Hence it is formally rigid.
For the classical current algebras $\gb=\g\otimes \C[z^{-1},z]$ 
with $\g$ a finite-dimensional simple Lie algebra,
Lecomte and Roger \cite{LeRog} showed that $\gb$ is formally
rigid.
Nevertheless, for both types of algebras, including their central
extensions, we obtained  deformations which are both locally 
geometrically and analytically non-trivial
\cite{FiaSchlwitt}, \cite{FiaSchlaff}. Hence they are not rigid in
the geometric and analytic sense. 
These families will be described in the following.



\section{Krichever-Novikov algebras}
\label{S:kn}

Our geometric families will be families of algebras of 
Krichever-Novikov type. These
algebras play an important  role in a global operator
approach to higher genus Conformal Field Theory.

They are  generalisations of the Virasoro
algebra, the current algebras and all their related algebras.
Let $M$ be a compact Riemann surface of genus $g$, or in terms
of algebraic geometry, a smooth projective curve over $\C$.
Let $N,K\in\N$ with $N\ge 2$ and $1\le K<N$. Fix
$$
I=(P_1,\ldots,P_K),\quad\text{and}\quad
O=(Q_1,\ldots,Q_{N-K})
$$
disjoint  ordered tuples of  distinct points (``marked points'',
``punctures'') on the
curve. In particular, we assume $P_i\ne Q_j$ for every
pair $(i,j)$. The points in $I$ are
called the {\it in-points}, the points in $O$ the {\it out-points}.
Sometimes we consider $I$ and $O$ simply as sets and set
$A=I\cup O$ as a set.

Here we will need
the following algebras.
Let $\A$ be the associative algebra 
consisting of those meromorphic functions 
on $M$ which are holomorphic outside the set of points $A$ with point-wise
multiplication.
Let $\L$ be the Lie algebra consisting of those meromorphic vector fields
which are holomorphic outside of $A$ with the usual Lie bracket of
vector fields.
The algebra $\L$ is called the {\it vector field algebra of 
Krichever-Novikov type}.
In the two point case they were introduced  by 
Krichever and Novikov \cite{KNFac}.
The corresponding generalisation to the multi-point case was done in
\cite{SchlDiss}, \cite{Schlkn}, \cite{Schleg}, \cite{Schlce}.
Obviously, both $\A$ and $\L$ are infinite dimensional algebras.

Furthermore, we will need the {\it higher-genus multi-point current
algebra of Krichever-Novikov type}.
We start with 
a complex finite-dimensional Lie algebra $\g$ and  endow
the tensor product $\Gb=\g\otimes_\C \A$ with the Lie bracket
\begin{equation}
[x\otimes f, y\otimes g]=[x,y]\otimes f\cdot g,
\qquad  x,y\in\g,\quad f,g\in\A.
\end{equation} 
The algebra  $\Gb$ is the higher genus current algebra.
It is an infinite dimensional Lie algebra and might be considered 
as the Lie algebra of $\g$-valued meromorphic functions on the
Riemann surface with only poles outside of $A$.

\medskip
The classical genus zero and $N=2$ point case is give by the 
geometric data
\begin{equation}\label{E:class}
M=\P^1(\C)=S^2, \quad I=\{z=0\},\quad 
 O=\{z=\infty\}.
\end{equation}
In this case the algebras are the well-known algebras of 
Conformal Field Theory (CFT).
For the function algebra
we obtain
$\A=\C[z^{-1},z]$, the algebra of Laurent polynomials.
The vector field algebra $\L$ is the Witt algebra
$\W$ generated by 
\begin{equation}
l_n=z^{n+1}\ddz,\ n\in\Z,
\quad
\text{ with Lie bracket}
\quad 
[l_n,l_m]=(m-n)\,l_{n+m}.
\end{equation}
The current algebra $\Gb$ is
the standard current algebra
$\gb=\g\otimes \C[z^{-1},z]$  with
Lie bracket
\begin{equation}
[x\otimes z^n, y\otimes z^m]=[x,y]\otimes z^{n+m},
\qquad  x,y\in\g,\quad n,m\in\Z.
\end{equation}

\medskip

In the classical situation the algebras are obviously graded by
taking as degree $\deg l_n:=n$ and $\deg x\otimes z^n:=n$.
For higher genus there is usually no grading.
But it was observed by Krichever and Novikov 
in the two-point case that a weaker
concept, an almost-graded structure, will be enough to develop an
interesting  theory of representations (Verma modules, etc.).
Let $\A$ be an (associative or Lie) algebra admitting a direct
decomposition as vector space $\ \A=\bigoplus_{n\in\Z} \A_n\ $.
The algebra $\A$ is called an {\it almost-graded}
algebra if (1) $\ \dim \A_n<\infty\ $ and (2)
there are constants $R$ and  $S$ such that 
\begin{equation}\label{E:eaga}
\A_n\cdot \A_m\quad \subseteq \bigoplus_{h=n+m+R}^{n+m+S} \A_h,
\qquad\forall n,m\in\Z\ .
\end{equation}
The elements of $\A_n$ are called {\it homogeneous  elements of degree $n$}.
By exhibiting a special basis,
for the multi-point situation such an almost grading was introduced in
\cite{SchlDiss,Schlkn, Schleg,Schlce}.
Essentially, this is done by fixing the order of the basis elements
at the points in $I$ in a certain manner and
in $O$ in a complementary way to make them unique
up to scaling.
In the following we will give an explicit description of the basis
elements
for those  genus zero and one situation we need.
Hence, we will not recall their general definition but only refer to
the above quoted articles.
\begin{proposition}
\cite{SchlDiss,Schlce}
The algebras $\L$, $\A$, and $\Gb$ 
are almost-graded.
The almost-grading depends on the splitting $A=I\cup O$.
\end{proposition}

\bigskip

In the construction of infinite dimensional representations of
these algebras with certain desired properties
(generated by a vacuum, irreducibility, unitarity, etc.)
one is typically forced to ``regularize'' a ``naive'' action
to make it well-defined.
Examples of importance in CFT are  
the fermionic Fock space representations which are
constructed by taking semi-infinite forms of a fixed weight.

From the mathematical point of view,
with the
help of a prescribed  procedure one modifies the action to 
make it well-defined.  On the other hand, one has to accept 
that the modified action 
in compensation will
only be  a projective Lie action. Such projective actions are
honest Lie actions  for suitable centrally extended algebras.
In the classical case they are well-known.
The unique non-trivial (up to equivalence and rescaling) central extension
of the Witt algebra $\W$ is the Virasoro algebra $\V$:
\begin{equation}\label{E:Vir}
[l_n,l_m] = (m-n)l_{n+m} + \tfrac{1}{12}(m^3-m)\delta_{n,-m}\,t,
\qquad
[l_n,  t] =0.
\end{equation}
Here $t$ is an additional element of the central extension which 
commutes with all other elements.
For the current algebra $\g\otimes \C[z^{-1},z]$ for $\g$ a simple Lie
algebra 
with Cartan-Killing form $\beta$, it is the corresponding affine Lie algebra 
$\gh$ (or, untwisted affine Kac-Moody algebra):
\begin{equation}\label{E:classcur}
[x\otimes z^n, y\otimes z^m]=[x,y]\otimes z^{n+m}-\beta(x,y)\cdot 
n\cdot\delta_{m}^{-n}\cdot t,
\quad [t,\gh]=0,
\quad x,y\in\g,\ n,m\in\Z.
\end{equation}
The additional terms in front of the elements $t$ are 2-cocycles of the
Lie algebras with values in the trivial module $\C$.
Indeed for a Lie algebra $V$ central extensions are 
classified (up to equivalence) by the second Lie algebra cohomology
$\H^2(V,\C)$
of $V$ with values in the trivial module $\C$.
Similar to the above, the  
bilinear form $\psi:V\times V\to\C$ is called  Lie algebra 2-cocycle
iff $\psi$ is antisymmetric and fulfills the cocycle condition
\begin{equation}
0=d_2\psi(x,y,z):=
\psi([x,y],z)+
\psi([y,z],x)+
\psi([z,x],y).
\end{equation}
It will be a coboundary if there exists a linear form  
$\kappa:V\to\C$ such that
\begin{equation}
\psi(x,y)=(d_1\kappa)(x,y):=\kappa([x,y]).
\end{equation}
To extend the classical cocycles to the Krichever-Novikov type
algebras they first have to be given in geometric terms.
Geometric  versions of the 2-cocycles 
are given as follows (see \cite{Schlloc} and \cite{Schlaff} for details).
For the vector field algebra $\L$ we take 
\begin{equation}\label{E:vecg}
\gamma_{S,R}(e,f):=\cintl{C_S} \left(\frac 12(e'''f-ef''')
-R\cdot(e'f-ef')\right)dz\ .
\end{equation}
Here the integration path $C_S$ is a 
loop separating the in-points from the out-points and
$R$ is a holomorphic projective connection (see
\cite[Def. 4.2]{FiaSchlwitt}) to make the integrand well-defined.
For the current algebra $\Gb$ we take
\begin{equation}\label{E:daff}
\gamma_{S}(x\otimes f,y\otimes g)=\beta(x,y)\cins fdg.
\end{equation}
The reader should be warned. 
For the classical algebras, i.e. the Witt and the current algebras
for the simple Lie algebras $\g$, there exists up to rescaling and equivalence
only one non-trivial central extension. This will be the Virasoro 
algebra for the Witt algebra and the affine Kac-Moody algebra for
the current algebra respectively.
This is not true anymore for higher genus or/and the multi-point
situation.
But it was shown in \cite{Schlloc} and \cite{Schlaff} that 
(again up to equivalence and rescaling) 
there exists only one non-trivial central extension 
which allows to extend the almost-grading by giving 
the element $t$ a degree in such a way that it will also be
almost-graded. This unique extension will be given by the
geometric cocycles
\refE{vecg}, \refE{daff}.


%
\section{The geometric families}\label{S:family}
\subsection{Complex torus}
Let $\tau\in\C$ with $\Im\tau>0$ and
$L$ be the lattice 
\begin{equation}
L=\langle 1,\tau\rangle_{\Z}:=\{m+n\cdot \tau\mid  m,n\in\Z\}\subset \C.
\end{equation} 
The  complex one-dimensional torus  is the quotient
$T=\C/L$. It carries a natural structure of a complex 
manifold coming from the structure of $\C$. It will be 
a compact  Riemann
surface of genus 1.

The field of meromorphic functions on $T$ is
generated by the doubly-periodic Weierstra\ss\ $\wp$ function and
its derivative $\wp'$ fulfilling the differential equation
\begin{equation}\label{E:diffequ}
(\wp')^2=4(\wp-e_1)(\wp-e_2)(\wp-e_3)=
4\wp^3-g_2\wp-g_3,
\end{equation}
with
\begin{equation}\label{E:delta}
\Delta:={g_2}^3-27{g_3}^2=16(e_1-e_2)^2
(e_1-e_3)^2(e_2-e_3)^2\ne 0.
\end{equation}
Furthermore,
\begin{equation}
g_2=-4(e_1e_2+e_1e_3+e_2e_3),\quad
g_3=4(e_1e_2e_3).
\end{equation}
The numbers $e_i$ are pairwise distinct, can be given as  
\begin{equation}
\wp(\frac {1}{2})=e_1,\qquad
\wp(\frac {\tau}{2})=e_2,\qquad
\wp(\frac {\tau+1}{2})=e_3,
\end{equation}
and fulfill
\begin{equation}
e_1+e_2+e_3=0.
\end{equation}

The function $\wp$ is an even meromorphic function 
with poles of order two at the
points of the lattice and holomorphic elsewhere.
The function $\wp'$ is an odd meromorphic function with poles of order 
three at the
points of the lattice  and holomorphic elsewhere.
$\wp'$ has zeros of order one at the points $1/2,\tau/2$ and $(1+\tau)/2$
and all their translates under the lattice.

We have to pass here to the algebraic-geometric picture.
The map
\begin{equation}
T\to\P^2(\C),\quad
z \bmod L\mapsto
\begin{cases}  
(\wp(z):\wp'(z):1),&z\notin L
\\
(0:1:0),&z\in L
\end{cases}
\end{equation}
realizes $T$ as a complex-algebraic smooth curve in the projective 
plane.
As its  genus is one it is  an elliptic curve. 
The affine coordinates are $X=\wp(z,\tau)$ and $Y=\wp'(z,\tau)$.
From \refE{diffequ} it follows that the affine part 
of the curve can be given by the
smooth cubic curve defined by
\begin{equation}\label{E:diffalg}
Y^2=4(X-e_1)(X-e_2)(X-e_3)=
4X^3-g_2X-g_3=:f(X).
\end{equation}
The point at infinity on the curve is the point
$\infty=(0:1:0)$.

We consider the algebras of Krichever-Novikov type corresponding to
the elliptic curve
and possible poles at $\bar z=\bar 0$ and $\bar z=\overline{1/2}$
\footnote{ Here $\bar z$ does not denote conjugation, but taking the
residue class modulo the lattice.}
(respectively 
in the algebraic-geometric picture, 
at the points $\infty$ and $(e_1,0)$).

\subsection{Vector field algebra}
First we consider the vector field algebra $\L$.
A basis of the vector field algebra 
is given by 
\begin{equation}\label{E:evec}
V_{2k+1}:=(X-e_1)^{k}Y\ddX,\quad
V_{2k}:=\frac 12f(X)(X-e_1)^{k-2}\ddX,\qquad k\in\Z.
\end{equation}
If we vary the points $e_1$ and $e_2$ 
(and accordingly $e_3=-(e_1+e_2)$) we obtain families
of curves and associated families of vector field algebras.
At least this is the case 
 as long as
the curves are non-singular.
To describe the families  in detail consider the following 
straight lines
\begin{equation}
D_s:=\{(e_1,e_2)\in\C^2\mid 
e_2=s\cdot e_1\},\   s\in \C,
\qquad
D_\infty:=\{(0,e_2)\in\C^2\}, 
\end{equation} 
and the open subset
\begin{equation} 
B=\C^2\setminus (D_1\cup D_{-1/2}\cup D_{-2})\subset\C^2.
\end{equation}
The curves are non-singular exactly over the points of $B$.
Over the exceptional $D_s$ at least two of the $e_i$ 
are the same.
For the vector field algebra we obtain
\begin{equation}\label{E:structell}
[V_n,V_m]=
\begin{cases}
(m-n)V_{n+m},&n,m \ \text{odd},
\\
(m-n)\big(V_{n+m}+3e_1V_{n+m-2}
\\
\qquad +\eee    V_{n+m-4}\big),&n,m \ \text{even},
\\
(m-n)V_{n+m}+(m-n-1)3e_1V_{n+m-2}
\\
\qquad +(m-n-2)\eee V_{n+m-4},&n\
\text{odd},\ m\ \text{even}.
\end{cases} 
\end{equation}
In fact these relations define Lie algebras for every pair $(e_1,e_2)\in\C^2$.
 We denote by
$\L^{(e_1,e_2)}$ the Lie algebra corresponding to $(e_1,e_2)$.
Obviously, $\L^{(0,0)}\cong \W$.
\begin{proposition} (\cite[Prop.~5.1]{FiaSchlwitt})
For $(e_1,e_2)\ne (0,0)$ the algebras $\L^{(e_1,e_2)}$ are not
isomorphic to the Witt algebra $\W$, but
$\L^{(0,0)}\cong \W$.
\end{proposition}
If we restrict our two-dimensional family to a line
$D_s$ ($s\ne\infty$) then we obtain a one-dimension family 
\begin{equation}\label{E:structs}
[V_n,V_m]=
\begin{cases}
(m-n)V_{n+m},&n,m \ \text{odd},
\\
(m-n)\big(V_{n+m}+3e_1V_{n+m-2}
\\
\qquad +e_1^2(1-s)(2+s)     V_{n+m-4}\big),&n,m \ \text{even},
\\
(m-n)V_{n+m}+(m-n-1)3e_1V_{n+m-2}
\\
\qquad
+(m-n-2)e_1^2(1-s)(2+s)    V_{n+m-4},&n\
\text{odd},\ m\ \text{even}.
\end{cases} 
\end{equation}
Here $s$ has a fixed value and  
 $e_1$ is the deformation parameter.
(A similar family exists for $s=\infty$.)
It can be shown that as long as $e_1\ne 0$ the algebras 
over two points in $D_s$ are
pairwise isomorphic but not isomorphic to the algebra
 over $0$, which is the Witt algebra.
Using the result $\H^2(\W,\W)=\{0\}$ of Fialowski \cite{Fia2}  we get
\begin{theorem}
Despite its infinitesimal and formal 
rigidity the Witt algebra $\W$ admits deformations $\L_t$ over the
affine line with $\L_0\cong \W$ which restricted to every 
(Zariski or analytic)
neighbourhood of $t=0$ are non-trivial.
\end{theorem}
The one-dimensional families \refE{structs} are examples of jump
deformations as $\L_t\cong\L_{t'}$ for $t, t'\ne 0$.
The isomorphism is given by rescaling 
the basis elements. This  is possible as long
as $e_1\ne 0$. In fact, using 
$V_n^*=(\sqrt{e_1})^{-n}V_n$  (for $s\ne \infty$) we obtain
for $e_1$ always the algebra with $e_1=1$ in the structure 
equations \refE{structs}.

Using the cocycle \refE{vecg} 
in the families \refE{structell}, \refE{structs}
a central term can be easily incorporated.
With respect to the flat coordinate $z-a$ we can take the projective
connection $R\equiv 0$.
The integral along a separating cocycle $C_S$ is obtained by 
taking the residue at $z=0$. In this way
we obtain geometric families of deformations for the Virasoro algebra.
They are locally non-trivial despite the fact that  the Virasoro
algebra is formally rigid.

\subsection{A family which is not a jump deformation}
Let us  stress the fact that the two-dimensional family
\refE{structell} is  not a jump deformation. 
In fact there exists even one-dimensional deformations as subfamilies
which are not jump deformations.
Take for example the smooth rational curve given by
\begin{equation}
C:=\{(e_1,e_2)\in \C\mid e_2=2e_1^2,\  e_1\in\C\}.
\end{equation} 
The rational parameter will  be $e_1$.
Automatically we have $e_3=-(1+2e_1)e_1$.

The curve passes through $(0,0)$. 
For every line $D_s$ there will be just one other point of intersection
with $C$. Its parameter value is given by $e_1=1/2 s$. Hence for $e_1$ small
the curve will not meet the exceptional lines $D_s,\ s=1,-1/2,-2$ a second
time. The curves corresponding to small $e_1\ne0$ values will
be nonsingular cubics, i.e. elliptic curves.

The elliptic modular function classifying elliptic curves up to
isomorphisms is given by 
\begin{equation}
j=1728\;\frac {g_2^3}{\Delta}.
\end{equation}
Expressing the $j$-function in terms of $e_i$, $i=1,2,3$ and 
substituting $e_3=-(e_1+e_2)$ yields
\begin{equation}
j(e_1,e_2)=1728\;
\frac {4(e_1^2+e_1e_2+e_2^2)^3}
{(e_1-e_2)^2(2e_1+e_2)^2(e_1+2e_2)^2},
\end{equation}
First we evaluate this along $D_s$, $s\ne 1,-1/2,-2$ and obtain
\begin{equation}
j(s,e_1)=1728\;\frac {4(1+s+s^2)^3}{(1-s)^2(2+s)^2(1+2s)^2},
\quad j(\infty)=1728.
\end{equation}
As this does not depend on the parameter $e_1$ anymore
it follows that the elliptic curves over $D_s\setminus\{0\}$ for
a fixed $s$ are isomorphic, in accordance with the fact that 
these deformations yield jump deformations.
A remark aside: the poles of $j(s)$ correspond exactly to the
exceptional lines. They correspond to nodal cubics, see \refS{inter}.

Next we evaluate  $j$ along the curve $C$ and obtain
\begin{equation}
j(e_1)=1728\;\frac {(1+2e_1+4e_1^2)^3}{(1-2e_1)^2(1+e_1)^2(1+4 e_1)^2}.
\end{equation}
This value will not be constant along $C$.
Furthermore, for small $e_1$ the values will be different.
This implies that the elliptic curves will be pairwise 
non-isomorphic.  
And the vector field algebras along this curve 
in the neighbourhood of $0$ will also be pairwise non-isomorphic.

By specializing 
$e_2=2e_1^2$ and $e_3=-e_1(1+2e_1)$ in \refE{structell}
we obtain
for the vector field algebra
\begin{equation}\label{E:structellc}
[V_n,V_m]=
\begin{cases}
(m-n)V_{n+m},&n,m \ \text{odd},
\\
(m-n)\big(V_{n+m}+3e_1V_{n+m-2}
\\
\qquad +  2e_1(1-2e_1)(1+e_1)  V_{n+m-4}\big),&n,m \ \text{even},
\\
(m-n)V_{n+m}+(m-n-1)3e_1V_{n+m-2}
\\
\qquad +(m-n-2)2e_1(1-2e_1)(1+e_1) V_{n+m-4},&n\
\text{odd},\ m\ \text{even}.
\end{cases} 
\end{equation}

\subsection{The current algebra}
Let $\g$ be a simple finite-dimensional Lie algebra (similar results 
are true for general Lie algebras) 
and $\A$   the algebra 
of meromorphic functions 
corresponding to the geometric situation discussed above.
A basis for  $\A$ is given by 

\begin{equation}
A_{2k}=(X-e_1)^{k},
\qquad
A_{2k+1}=\frac 12 Y\cdot (X-e_1)^{k-1}, \qquad k\in\Z.
\end{equation}

We calculate for the elements of $\Gb$
\begin{equation}\label{E:structaff}
[x\otimes A_n, y\otimes A_m]
=\begin{cases}
[x,y]\otimes A_{n+m},& \text{$n$ or $m$ even},
\\
[x,y]\otimes A_{n+m}+3e_1 [x,y]\otimes A_{n+m-2}
\\
\ +(e_1-e_2)(2e_1+e_2) [x,y]\otimes A_{n+m-4},
& \text{$n$ and $m$ odd}.
\end{cases}
\end{equation}
If we let $e_1$ and $e_2$ (and hence also $e_3$) go to zero
we obtain the classical current algebra as degeneration.
Again it can be shown that the family, even if restricted on $D_s$, 
is locally  non-trivial, see \cite{FiaSchlaff}. 
Recall that by results of Lecomte and Roger \cite{LeRog} 
the current algebra is formally rigid if $\g$ is simple.
But our families show that it is neither geometrically nor
analytically rigid.

The families over $D_s$ are jump deformations. But again there exists
one-dimensional subfamilies of deformations of \refE{structaff} which
are non-trivial and not jump deformations. 

Also in this case we can construct families of centrally extended algebras by
considering the cocycle 
\refE{daff}. In this way we obtain non-trivial 
deformation families for the formally rigid classical affine
algebras of Kac-Moody type.
The cocycle \refE{daff} is 
\begin{equation}
\label{E:famcent}
\gamma(x\otimes A_n,y\otimes A_m)=
p(e_1,e_2)\cdot \beta(x,y)\cdot\cins A_ndA_m.
\end{equation}
Here $p(e_1,e_2)$ is an arbitrary polynomial in the variables
$e_1$ and $e_2$.
 and $\beta$ the Cartan-Killing form.
The integral can be calculated \cite[Thm. 4.6]{FiaSchlaff}
as
\begin{equation}
\label{E:famcent1}
\cins A_ndA_m=
\begin{cases}
\qquad -n\delta_m^{-n},&\text{$n$, $m$ even},
\\
\qquad\qquad 0,&\text{$n$, $m$ different parity},
\\
-n\delta_m^{-n}+
3e_1(-n+1)\delta_m^{-n+2}+
\\
+(e_1-e_2)(2e_1+e_2)(-n+2)\delta_m^{-n+4}
,&\text{$n$, $m$ odd}.
\end{cases}
\end{equation}


\section{The geometric interpretation}
\label{S:inter}
If we take $e_1=e_2=e_3$ in the definition of the
cubic curve \refE{diffalg} we obtain the cuspidal cubic $E_C$ 
with affine part given by the polynomial 
$Y^2=4X^3$. It has a singularity at $(0,0)$ and the
desingularization is given by the projective line $\P^1(\C)$.
This says there exists a surjective (algebraic) map 
$\pi_C:\P^1(\C)\to E_C$ which outside the singular point is
$1:1$. Over the cusp lies
exactly one point.
The vector fields, resp. the functions, resp. the 
$\g$-valued functions can be degenerated to 
$E_C$ and pull-backed to vector fields, resp. functions, resp.
$\g$-valued functions  on 
$\P^1(\C)$.
The  point $(e_1,0)$ where a pole is allowed moves to
the cusp. The other point stays at infinity.
In particular by pulling back the degenerated vector 
field algebra we obtain
the algebra of vector fields with two possible poles,
which is the Witt algebra.
And by pulling back the degenerated current algebra we obtain 
the classical current algebra.

\medskip

The exceptional lines $D_s$ for $s=1,-1/2,-2$ are related
to interesting geometric situations.
Above $D_s\setminus\{(0,0)\}$ with these  values of $s$, two of
the $e_i$ are the same, the third one remains distinct.
The curve will be a nodal cubic $E_N$ defined by
$Y^2=4(X-e)^2(X-e)$. The singularity will be a node with the
coordinates $(e,0)$.
Again the desingularization will be the projective line
$\pi_N:\P^1(\C)\to E_N$. But now above the node there will be
two points in $\P^1(\C)$.
For the pull-backs  we have the following two
situations:

{\bf (1)} If $s=1$ or $s=-2$ then $e=e_1$ and
 the node is a possible point for a pole.
We obtain objects  on $\P^1(\C)$ which  beside
the pole at $\infty$ might have poles at two other places. Hence, we
obtain a three-point Krichever-Novikov  algebras of genus 
0.

{\bf (2)} If $s=-1/2$ then  at the node there is no pole.
The number of possible poles for the pull-back remains two.
We obtain certain subalgebras of the classical two point case.
Additionally, for the vector field case we have to pay attention
to the fact that the vector fields obtained 
by pull-back acquire zeros at the
points lying above the node.

These algebras were identified and studied in detail in 
\cite{FiaSchlwitt}, \cite{FiaSchlaff}  and \cite{Schltori}.

The deformed families are of importance 
for the quantization of conformal field theories. In this context 
the behaviour of objects when 
we approach 
the boundary of the moduli space of curves with marked 
points has to be studied.

\appendix
\section{Definition of a deformation}\label{S:defin}

In the following we will assume that $A$ is a
commutative algebra over  $\K$ (where $\K$ is
a field  of  characteristic zero)
which  admits an augmentation 
 $\epsilon: A \to \K$.
This says that $\epsilon$ is a $\K$-algebra homomorphism, e.g.
$\epsilon(1_A)=1$.
The ideal $m_\epsilon:= \text{Ker\,} \epsilon$ is 
a maximal ideal of $A$. 
Vice versa, given a maximal ideal $m$ of $A$ with 
$A/m\cong\K$, the natural quotient map defines an
augmentation. 

If $A$ is a finitely generated $\K$-algebra over
an algebraically closed field $\K$ then $A/m\cong\K$
is true for every maximal ideal $m$. Hence, in this case 
every such $A$ admits at least one augmentation and
all maximal ideals are coming from augmentations.

Let us consider a Lie algebra $\L$ over the field $\K$, 
 $\epsilon$ a fixed augmentation of $A$, 
and
$m=\text{Ker\,} \epsilon$ the associated maximal ideal. 
\begin{definition}
\label{D:glob} (\cite{FiaFuMini})
 A \emph{deformation} $\lambda$ of $\L$ with 
base $(A,m)$ or simply with  base $A$, is a Lie $A$-algebra structure
on the tensor product $A\tensor_\K\L$ with bracket $[.,.]_\lambda$ such
that
\begin{equation}
\epsilon \textstyle{\tensor} \text{id} : A \textstyle{\tensor} \L \to
 \K \textstyle{\tensor} \L = \L
\end{equation}
is a Lie algebra homomorphism.

Specifically, it means that for all $a,b \in A$ and $x,y \in \L$,
\smallskip

\quad (1) $[a \tensor x, b \tensor y]_\lambda = (ab \tensor \text{id\,})[1 \tensor x,
1 \tensor y]_\lambda$,
\smallskip

\quad (2) $[.,.]_\lambda$ is skew-symmetric and satisfies the Jacobi 
identity,
\smallskip

\quad (3) $\epsilon \tensor \text{id\,}([1 \tensor x, 1 \tensor y]_\lambda) = 1
\tensor [x,y]$.
\end{definition}

\noindent 
By Condition (1) to describe a deformation 
it is enough to give  the elements 
\newline
$[1 \tensor x, 1
\tensor y]_\lambda$ for all $x,y \in \L$. 
If $B=\{z_i\}_{i\in J}$ is a basis of $\L$ it follows 
from Condition (3) that the Lie product has the form
\begin{equation}\label{E:defex}
[1 \textstyle{\tensor} x, 1 \textstyle{\tensor} y]_\lambda = 1
\textstyle{\tensor} [x,y] + \sum_i' a_i \textstyle{\tensor} z_i,
\end{equation}
with $a_i=a_i(x,y) \in m$, $z_i \in B$.
Here $\sum'$ denotes a finite sum.  
Clearly, Condition (2) is an additional condition which has to
be satisfied.

If we use $A=\C[t]$ we get exactly the notion of a one parameter geometric
deformation discussed above.

A deformation is called \emph{trivial} if $A\tensor_{\K} \L$
carries the trivially extended Lie structure, i.e. \refE{defex} reads as 
$[1\tensor x,1\tensor y]_\lambda=1\tensor[x,y]$.

Two deformations of a Lie algebra $\L$ with the same base $A$ are
called \emph{equivalent} if there exists a Lie algebra isomorphism
between the two copies of $A \tensor \L$ with the two Lie algebra
structures, compatible with $\epsilon \tensor \text{id}$.

\medskip

Formal deformations are defined in a similar way.
Let $A$ be a complete local algebra over $\K$, so
$A = \overleftarrow{\lim\limits_{n\to \infty}}(A/m^n)$, where $m$ is
the maximal ideal of $A$. Furthermore,  we 
assume that $A/m\cong\K$, and $\dim (m^k/m^{k+1})<\infty$ for all $k$.
\begin{definition} (\cite{Fia1})
A \emph{formal deformation} of $\L$ with base $A$ is a Lie algebra
structure on the completed tensor product
$A \htensor\L =
    \overleftarrow{\lim\limits_{n\to \infty}}((A/m^n)\tensor \L)$
such that
\begin{equation}
\epsilon\textstyle{\htensor}
\text{id}: A \textstyle{\htensor}\L \to \K\tensor \L = \L
\end{equation}
is a Lie algebra homomorphism.
\end{definition}
If $A=\C[[t]]$, then a formal deformation of $\L$ with base $A$
is the same as a formal one parameter deformation discussed 
above.
There is an analogous definition for equivalence of deformations
parameterized by a complete local algebra.



\end{document}